\documentclass[a4paper,11pt]{article} 
\usepackage{float}
\usepackage{booktabs} 
\usepackage{tikz-qtree}
\usetikzlibrary{arrows.meta}
\usepackage[a4paper, margin=.9in]{geometry}
\usepackage{amsmath,amsthm,amsfonts,amssymb,array}
\usepackage[square,sort,comma,numbers]{natbib}
\usepackage{subeqnarray}
\usepackage{blkarray}
\usepackage[utf8]{inputenc}
\usepackage[english]{babel}
\usepackage{mathabx}
\usepackage{mathrsfs}
\usepackage{colortbl}
\usepackage{graphicx,epsfig}
\newcounter{myeqno}
\usepackage{color}
\usepackage{tcolorbox}
\usepackage{ upgreek }
\usepackage{relsize}
\usepackage{multirow}
\usepackage{MnSymbol}
\definecolor{shadecolor}{gray}{0.75}

\theoremstyle{definition}

\newtheorem{theorem}{Theorem}[section]
\newtheorem{example}{Example}[section]
\newtheorem{corollary}{Corollary}[section]

\newtheorem{lemma}{Lemma}[section]

\usepackage[misc,geometry]{ifsym}
\date{}
\usepackage[colorlinks=true]{hyperref}
\hypersetup{urlcolor=blue, citecolor=blue}
\newlength{\defbaselineskip}
\newcommand{\setlinespacing}[1]%
{\setlength{\baselineskip}{#1 \defbaselineskip}}

\begin{document}

		\title{\textbf{Resistance distance and Kirchhoff index in central vertex join and central edge join of two graphs}}

\author{Haritha T$^1$\footnote{harithathottungal96@gmail.com},  Chithra A. V$^1$\footnote{chithra@nitc.ac.in}
	\\ \\ \small 
 1 Department of Mathematics, National Institute of Technology Calicut,\\\small Calicut-673 601, Kerala, India\\ \small}

\maketitle	
\begin{abstract}
The central graph $C(G)$ of a graph $G$ is the graph obtained by inserting a new vertex into each edge of $G$ exactly once and joining all the non-adjacent vertices in $G$. Let $G_1$ and $G_2$ be two vertex disjoint graphs. The central vertex join of $G_1$ and $G_2$ is the graph  $ G_1\dot{\vee} G_2$, is obtained from $C(G_1)$ and $G_2$ by joining each vertex of $G_1$ with every vertex of $G_2$. The central edge join of $G_1$ and $G_2$ is the graph  $ G_1\veebar  G_2$, is obtained from $C(G_1)$ and $G_2$ by joining each vertex corresponding to the edges of $G_1$ with every vertex of $G_2$. In this article, we obtain formulae for the resistance distance and Kirchhoff index of $G_1\dot{\vee} G_2$ and $ G_1\veebar  G_2$. In addition, we provide the resistance distance, Kirchhoff index, and Kemeny's constant of the central graph of a graph.
\end{abstract}
\section{Introduction}
Throughout this article, we consider simple connected undirected graphs. Let $G$ be a graph consisting of $n$ vertices $V(G)= \{v_1, v_2, \ldots, v_n\}$ and $m$ edges $E(G)= \{e_1, e_2, \ldots, e_m\}$. The \textit{adjacency matrix} $A(G)$ of a graph $G$ is an $n\times n$ symmetric matrix defined as $$A(G)= (a_{ij})_{n\times n}= \begin{cases}1, &\text{if $v_i$ is adjacent to $v_j$,}\\0, &\text{otherwise}.\end{cases}$$  Let the all-one entry matrix be denoted by $J$ in appropriate order, and the identity matrix by $I$. Denote the degree of a vertex $v_i$ in $G$ by $d_i$ and a column vector with all entries by $e$. Let $D(G)$ be the diagonal matrix of vertex degrees in $G$. Then the \textit{Laplacian matrix} of $G$ is defined by $L(G)=(l_{ij})= D(G)-A(G).$
Klein and Randić introduced the notion of \textit{resistance distance} in 1993 \cite{klein1993resistance}. The authors presented a new point of view, if each edge of a connected graph is assigned a fixed unit resistance, then the effective resistance between pairs of vertices is a graphical distance. For an $m\times n$ matrix $M$, the \textit{Moore-Penrose inverse} \cite{ben2003generalized} of $M$, denoted by $M^\dagger$, is the unique matrix 
$X$ such that $MXM= M$, $XMX= X$, $(MX)^T= MX$ and $(XM)^T= XM$.
     For an $m\times n$ matrix $M$, the matrix $P$ of order $n\times m$ is said to be a \textit{$\{ 1\}$-inverse} of $M$ (denoted by $M^{(1)}$) \cite{ben2003generalized} if $MPM = M$.

The Moore-Penrose inverse and the $\{1\}$-inverse of $L=(l_{ij})$ of the underlying graph $G$ are used to calculate the resistance distance $r_{ij}$ \cite{bapat2010graphs} between two vertices $v_i$ and $v_j$. The equation is as follows

 \begin{equation*}
     \begin{aligned}
         r_{ij} &= l^{(1)}_{ii}+ l^{(1)}_{jj}-l^{(1)}_{ij}-l^{(1)}_{ji}=l_{ii}^{\dagger}+l_{jj}^{\dagger}-2l_{ij}^{\dagger}.
     \end{aligned}
 \end{equation*}

The matrix $R(G)= (r_{ij})_{n\times n}$ is called the resistance matrix of $G$.

 The \textit{Kirchhoff index} of $G$, also known as the total resistance of a network, represented as $\mathcal{K}f(G)$ \cite{klein1993resistance,bonchev1994molecular},
is defined as,
$$\mathcal{K}f(G)= \sum_{i<j} r_{ij}.$$
\par \textit{Kemeny's constant}  \cite{kirkland2012group} is essential in the theory of random walks, and it measures the average time it takes for a random walk to reach a vertex. It is defined as
 $$\kappa(G)= \frac{1}{4m}\sum_{v_i, v_j\in V(G)}d_{i}d_{j}r_{ij}.$$

\par The importance of resistance distance extends to combinatorial matrix theory  \cite{bapat2010graphs, bapat2010resistance} as well as spectral graph theory \cite{bapat1999resistance,chen2007resistance,bapat2003simple,zhang2007resistance,ressub}. See \cite{evans2022algorithmic} for an overview of techniques for determining resistance distance in graphs.
The Kirchhoff index and resistance distance formulas for a variety of graph classes and graph operations were provided in \cite{ressub,liu2015resistance,wang2019resistance}. In order to investigate complicated networks with attributes abstracted from the real world, graph operations have been used extensively. The \textit{central graph} \cite{vivin2008harmonious} of $G$, denoted by $C(G)$, is obtained by inserting a new vertex into each edge of $G$ and joining all the non adjacent vertices in $G$. Using the central graph, Jahfar et al. defined two new graph operations called the central vertex join and central edge join of two graphs \cite{jahfar2020central}.
Let $G_1$ and $G_2$ be any two graphs, then the \textit{central vertex join} of $G_1$ and $G_2$ is the graph  $ G_1\dot{\vee} G_2$ obtained from $C(G_1)$ and $G_2$ by joining each vertex of $G_1$ with every vertex of $G_2$. The \textit{central edge join} of $G_1$ and $G_2$ is the graph  $ G_1\veebar  G_2$ obtained from $C(G_1)$ and $G_2$ by joining each vertex corresponding to the edges of $G_1$ with every vertex of $G_2$. The objective of this study is to provide the resistance distance and Kirchhoff index within the central vertex join and central edge join of two graphs.

\section{Preliminaries}
Through this section we present some useful lemmas and theorems.

\begin{lemma}{\cite{1inv}}\label{pinv}
    Let $L= \begin{bmatrix}
        L_1& L_2\\L_2^T& L_3
    \end{bmatrix}$ be the Laplacian matrix of a connected graph $G.$
If $L_3$ is non-singular, then
$$L^{(1)}= \begin{bmatrix}
    M^{\dagger}& -M^{\dagger}L_2L_3^{-1}\\
    -L_3^{-1}L_2^TM^{\dagger}& L_3^{-1}+L_3^{-1}L_2^TM^{\dagger}L_2L_3^{-1}
\end{bmatrix},$$
where $M= L_1-L_2L_3^{-1}L_2^T.$

\end{lemma}
\begin{lemma}\cite{ressub}\label{uv}
    Let $G$ be a connected graph. For any $v_i, v_j\in V(G)$,
    $$r_{ij}(G)= \frac{1}{d_i}\Bigl(1+\sum_{v_k\in N(v_i)}r_{kj}(G)-\frac{1}{d_i}\sum_{v_k, v_l\in N(v_i)}r_{kl}(G)\Bigr).$$

\end{lemma}
\begin{lemma}\cite{sun2015some}\label{Kf}
     Let $G$ be a connected graph of order $n$. Then
 $$Kf(G) = ntr\left(L^{(1)}(G)\right)- e^T L^{(1)}(G)e= ntr L^{\dagger}(G).$$
\end{lemma}
\begin{lemma}\cite{klein2002resistance}\label{rssum}
    For a connected graph $G$ of order $n$,
    $$\sum_{i<j, v_iv_j\in E(G)}r_{ij}(G)= n-1.$$
\end{lemma}

\section{Main results}
This section provides the resistance distance in central vertex(edge) join of two graphs and discusses its Kirchhoff index.\\

\begin{lemma}\label{aL}
    Let $L$ be the Laplacian matrix of a graph $G$ of order $n$. For any nonzero $a,b$, we have
    $$\left(aL+bI-\frac{b}{n}J\right)^{\dagger}= (aL+bI)^{-1}-\frac{1}{bn}J.$$
\end{lemma}
\begin{proof}
    Let $M= (aL+bI)^{-1}-\frac{1}{bn}J$. Then we have,\\
    \begin{equation*}
        \begin{aligned}
            M\left(aL+bI-\frac{b}{n}J\right)&= \left(aL+bI-\frac{b}{n}J\right)M= I-\frac{1}{n}J\\
            M\left(aL+bI-\frac{b}{n}J\right)M&= M\\
            \left(aL+bI-\frac{b}{n}J\right)M\left(aL+bI-\frac{b}{n}J\right)&= \left(aL+bI-\frac{b}{n}J\right).
            \end{aligned}
    \end{equation*}
    Therefore, $M= \left(aL+bI-\frac{b}{n}J\right)^{\dagger}.$
\end{proof}
\begin{theorem}\label{cent}
    For a connected graph $G$ of order $n$, the resistance distance between pairs of vertices in its central graph $C(G)$ is as follows,
    \begin{itemize}
        \item [(i)] for any two vertices $v_i, v_j\in V(G)$,
        $$r_{ij}(C(G))= \left(\frac{-1}{2}L(G)+nI\right)^{-1}_{ii}+\left(\frac{-1}{2}L(G)+nI\right)^{-1}_{jj}-2\left(\frac{-1}{2}L(G)+nI\right)^{-1}_{ij},$$
        \item[(ii)] for $v_i= v_{e_i}\in V(C(G))$, where $e_i= v_pv_q\in E(G), v_j\in V(G)$,
        $$r_{ij}(C(G))= \frac{1}{2}\left(1+r_{pj}(C(G))+r_{qj}(C(G))-\frac{1}{2}r_{pq}(C(G))\right),$$
        \item[(iii)] for $v_i= v_{e_i}, v_j= v_{e_j}\in V(C(G))$, where $e_i= v_pv_q, e_j= v_sv_t\in E(G),$
       $$r_{ij}(C(G))= 1+\frac{1}{4}\Bigl(r_{ps}(C(G))+r_{pt}(C(G))+r_{qs}(C(G))+r_{qt}(C(G))-r_{pq}(C(G))-r_{st}(C(G))\Bigr).$$
        
    \end{itemize}
\end{theorem}

\begin{proof}
    By a proper labelling of vertices in $C(G)$, its Laplacian matrix can be written as
$$L(C(G))= \begin{bmatrix}
    nI-J+A(G)& -Q(G)\\
    -Q(G)^T& 2I
\end{bmatrix}.$$
Let $L_1=  nI-J+A(G), L_2= -Q(G),$ and $L_3= 2I.$
Let \begin{equation*}
    \begin{aligned}
        M&= L_1-L_2L_3^{-1}L_2^T\\
        &= nI-J+A-\frac{1}{2}Q(G)Q(G)^T\\
        &= nI-J-\frac{1}{2}L(G).
    \end{aligned}
\end{equation*}
Now by Lemma \ref{aL}, $M^{\dagger}= (-\frac{1}{2}L(G)+nI)^{-1}-\frac{1}{n^2}J.$ Then by Lemma \ref{pinv} the $\{1\}$-inverse of
$L(C(G))$ is\\

    $$\begin{bmatrix}
     (-\frac{1}{2}L(G)+nI)^{-1}-\frac{1}{n^2}J& \frac{1}{2}(\frac{-1}{2}L(G)+nI)^{-1}Q(G)-\frac{1}{n^2}J\\&\\
     \frac{1}{2}Q(G)^T(\frac{-1}{2}L(G)+nI)^{-1}-\frac{1}{n^2}J&\frac{1}{2}I+\frac{1}{4}Q(G)^T(\frac{-1}{2}L(G)+nI)^{-1}Q(G)-\frac{1}{n^2}J
    \end{bmatrix}.$$\\
\noindent
 For any two vertices $v_i,\; v_j\in V(G)$,\\
        $$r_{ij}(C(G))= \left(\frac{-1}{2}L(G)+nI\right)^{-1}_{ii}+\left(\frac{-1}{2}L(G)+nI\right)^{-1}_{jj}-2\left(\frac{-1}{2}L(G)+nI\right)^{-1}_{ij}.$$\\

\noindent From Lemma \ref{uv}, for $v_i= v_{e_{i}}\in V(C(G))$, where $e_i= v_pv_q\in E(G_1)$, $v_j\in V(G),$\\

\begin{equation*}
    \begin{aligned}
     r_{ij}(C(G))&=  \frac{1}{2}\left(1+r_{pj}(C(G))+r_{qj}(C(G))-\frac{1}{2}r_{pq}(C(G))\right).   
    \end{aligned}
\end{equation*}
 \\
From Lemma \ref{uv}, for $v_i= v_{e_{i}}, v_j= v_{e_{j}}\in V(C(G)),$ where $e_i= v_pv_q, e_j= v_sv_t\in E(G),$\\
 \begin{equation*}
     \begin{aligned}
        r_{ij}(C(G))&= 1+\frac{1}{4}\Bigl(r_{ps}(C(G))+r_{pt}(C(G))+r_{qs}(C(G))+r_{qt}(C(G))-r_{pq}(C(G))-r_{st}(C(G))\Bigr).
     \end{aligned}
 \end{equation*}
        
\end{proof}
Next theorem gives a formula for the Kirchhoff index of the central graph of a graph.
\begin{theorem}\label{kfcg}
    The Kirchhoff index of the central graph of a connected graph G of order $n$ and size $m$ is
     \begin{equation*}
     \begin{aligned}
       Kf(C(G))&= (m+n)\left(tr \left(\frac{-1}{2}L(G)+nI\right)^{-1}\right)-\frac{1}{4}\alpha^T tr \left(\frac{-1}{2}L(G)+nI\right)^{-1}\alpha+\frac{m(m+n-1)}{2}-\frac{2m+n}{n}\\&\;\;\;+\frac{(m+n)(n-1)}{4}+\frac{m+n}{4}\Bigl(\sum_{i<j,v_iv_j\in E(G)}4\left(\frac{-1}{2}L(G)+nI\right)^{-1}_{ij}\Bigr),
     \end{aligned}
 \end{equation*}
    
    where $\alpha= (d_1, d_2, \ldots, d_n)^T.$
\end{theorem}
\begin{proof}
    By using Lemma \ref{Kf} and Theorem \ref{cent}, we have
    \begin{equation*}
        \begin{aligned}
            Kf(C(G))&= (m+n)\left(tr \left(\frac{-1}{2}L(G)+nI\right)^{-1}-\frac{1}{n}+\frac{m}{2}-\frac{m}{n^2}+\frac{1}{4}tr \left(Q^Ttr \left(\frac{-1}{2}L(G)+nI\right)^{-1}Q\right)\right)\\&\;\;\;-\frac{m}{2}+\frac{m^2}{n^2}-\frac{1}{4}e^TQ^Ttr \left(\frac{-1}{2}L(G)+nI\right)^{-1}Qe\\
            &= (m+n)\left(tr \left(\frac{-1}{2}L(G)+nI\right)^{-1}\right)-\frac{1}{4}\alpha^T tr \left(\frac{-1}{2}L(G)+nI\right)^{-1}\alpha+\frac{m(m+n-1)}{2}-\frac{2m+n}{n}\\&\;\;\;+\frac{m+n}{4}\left(\sum_{i<j,v_iv_j\in E(G)}\left(\frac{-1}{2}L(G)+nI\right)^{-1}_{ii}+\left(\frac{-1}{2}L(G)+nI\right)^{-1}_{jj}+2\left(\frac{-1}{2}L(G)+nI\right)^{-1}_{ij}\right)\\
            &= (m+n)\left(tr \left(\frac{-1}{2}L(G)+nI\right)^{-1}\right)-\frac{1}{4}\alpha^T tr \left(\frac{-1}{2}L(G)+nI\right)^{-1}\alpha+\frac{m(m+n-1)}{2}-\frac{2m+n}{n}\\&\;\;\;+\frac{m+n}{4}\Bigl(\sum_{i<j,v_iv_j\in E(G)}r_{ij}(G)+4\left(\frac{-1}{2}L(G)+nI\right)^{-1}_{ij}\Bigr).\\
        \end{aligned}
    \end{equation*}
    Now apply Lemma \ref{rssum} to get 
     \begin{equation*}
     \begin{aligned}
       Kf(C(G))&= (m+n)\left(tr \left(\frac{-1}{2}L(G)+nI\right)^{-1}\right)-\frac{1}{4}\alpha^T tr \left(\frac{-1}{2}L(G)+nI\right)^{-1}\alpha+\frac{m(m+n-1)}{2}-\frac{2m+n}{n}\\&\;\;\;+\frac{(m+n)(n-1)}{4}+\frac{m+n}{4}\Bigl(\sum_{i<j,v_iv_j\in E(G)}4\left(\frac{-1}{2}L(G)+nI\right)^{-1}_{ij}\Bigr),
     \end{aligned}
 \end{equation*}
    
    where $\alpha= (d_1, d_2, \ldots, d_n)^T.$
\end{proof}
\begin{corollary}
    For a connected graph $G$ of order $n$ and size $m$, let $\Delta$ denote the maximum degree of $G$. Then
    $$Kf(C(G))\leq (m+n) \frac{(\Delta+2)}{2} tr \left(\frac{-1}{2}L(G)+nI\right)^{-1}+\frac{m(2m+n+1)+8}{2}-\frac{n(n+1)}{2}-\frac{4m(m+2n)}{2n^2}.$$
\end{corollary}
\begin{corollary}
Let $G$ be a $d$-regular graph of order $n$, then
\begin{equation*}
\begin{aligned}
Kf(C(G))&= \frac{n(d+2)}{2}tr\left(\frac{-1}{2}L(G)+nI\right)^{-1}+\frac{n(d+2)}{8}tr \left(Q^T\left(\frac{-1}{2}L(G)+nI\right)^{-1}Q\right)-\frac{d(d+n)}{4}\\&\;\;\;+\frac{n^2d(d+2)}{8}-(d+1).
\end{aligned}
\end{equation*}
\end{corollary}
Next corollary gives the Kemeny's constant of central graph of a graph.
\begin{corollary}
    For a connected graph $G$ of order $n$ and size $m$,
    \begin{equation*}
        \begin{aligned}
            \kappa(C(G))&= \frac{1}{4m+2n(n-1)}\Bigl(\sum_{v_i,v_j\in V(G)}4\times r_{ij}(C(G))+\sum_{v_{e_i}\in V(C(G)),v_j\in V(G)}2(n-1)\times r_{ij}(C(G))\\&\;\;\;\;+\sum_{v_{e_i},v_{e_j}\in V(C(G))}(n-1)^2\times r_{ij}(C(G))\Bigr).
        \end{aligned}
    \end{equation*}
    
\end{corollary}

The next example provides the Kirchhoff index of central graph of some standard graphs.

\begin{example}
Using Theorem \ref{kfcg} we get the following:
    \begin{itemize}
        \item [(i)] Cycle $C_n$ ($n>2$)
        $$Kf(C(C_n))= 4n\sum_{k=1}^{n}\frac{1}{(n-1)+\operatorname{cos}\frac{2\pi k}{n}}+n^2-3.$$
        \item[(ii)] Regular complete bipartite graph $K_{n,n}$ 
        $$Kf(C(K_{n,n}))= \frac{120n^4+168n^3-157n^3-164n+488}{120}.$$
    \end{itemize}
\end{example}

\begin{theorem}\label{cvj}
    For $i= 1,2$, let $G_i$ be a graph of order $n_i$ and size $m_i$. Then
    \begin{itemize}
        \item [(i)] for any two vertices $v_i,\; v_j\in V(G_1)$,\\
        \begin{equation*}
            \begin{aligned}
                r_{ij}(G_1\dot{\vee} G_2)&= \left(-\frac{1}{2}L(G_1)+(n_1+n_2)I\right)^{-1}_{ii}+\left(-\frac{1}{2}L(G_1)+(n_1+n_2)I\right)^{-1}_{jj}\\
                &-2\left(-\frac{1}{2}L(G_1)+(n_1+n_2)I\right)^{-1}_{ij},
            \end{aligned}
        \end{equation*}
        
        \item[(ii)] for any two vertices $v_i,\; v_j\in V(G_2)$,\\
        $$r_{ij}(G_1\dot{\vee} G_2)= (n_{1}I+L(G_2))^{-1}_{ii}+(n_{1}I+L(G_2))^{-1}_{jj}-2(n_{1}I+L(G_2))^{-1}_{ij},$$
        \item[(iii)] for any $v_i\in V(G_1),\; v_j\in V(G_2),$\\
         $$r_{ij}(G_1\dot{\vee} G_2)= \left(-\frac{1}{2}L(G_1)+(n_1+n_2)I\right)^{-1}_{ii}+(n_{1}I+L(G_2))^{-1}_{jj}-\frac{1}{n_1(n_1+n_2)},$$

        \item[(iv)] for $v_i= v_{e_{i}}\in V(C(G_1))$, where $e_i= v_pv_q\in E(G_1)$, $v_j\in V(G_1)\cup V(G_2),$\\
        $$r_{ij}(G_1\dot{\vee} G_2)= \frac{1}{2}\Bigl(1+r_{pj}(G_1\dot{\vee} G_2)+r_{qj}(G_1\dot{\vee} G_2)-\frac{1}{2}r_{pq}(G_1\dot{\vee} G_2)\Bigr),$$
        \item[(v)] for $v_i= v_{e_{i}}, v_j= v_{e_{j}}\in V(C(G_1)),$ where $e_i= v_pv_q, e_j= v_sv_t\in E(G_1),$\\
        \begin{equation*}
            \begin{aligned}
                r_{ij}&= 1+\frac{1}{4}\Bigl(r_{ps}(G_1\dot{\vee} G_2)+r_{pt}(G_1\dot{\vee} G_2)+r_{qs}(G_1\dot{\vee} G_2)+ r_{qt}(G_1\dot{\vee} G_2)-r_{st}(G_1\dot{\vee} G_2)-r_{pq}(G_1\dot{\vee} G_2)\Bigr).
            \end{aligned}
        \end{equation*}

    \end{itemize}
\end{theorem}

\begin{proof}
By a proper labelling of vertices in $G_1\dot{\vee} G_2$, its Laplacian matrix can be written as
$$L(G_1\dot{\vee} G_2)= \begin{bmatrix}
    (n_1+n_2)I-J+A(G_1)& -Q(G_1)& -J\\
    -Q(G_1)^T& 2I& 0\\
    -J& 0& n_1I+L(G_2)
\end{bmatrix}.$$
Let $L_1= (n_1+n_2)I-J+A(G_1), L_2= \begin{bmatrix}
    -Q(G_1)& -J
\end{bmatrix},$ and $L_3= \begin{bmatrix}
    2I & 0\\0 & n_1I+L(G_2)
\end{bmatrix}.$
Let \begin{equation*}
    \begin{aligned}
        M&= L_1-L_2L_3^{-1}L_2^T\\
        &= (n_1+n_2)I-J+A(G_1)-\begin{bmatrix}
    -Q(G_1)& -J
\end{bmatrix}\begin{bmatrix}
    2I & 0\\0 & n_1I+L(G_2)
\end{bmatrix}\begin{bmatrix}
    -Q(G_1)^T\\-J^T
\end{bmatrix}\\&= (n_1+n_2)I-J+A(G_1)-\frac{1}{2}Q(G_1)Q(G_1)^T-\frac{n_2}{n_1}J\\
&= -\frac{1}{2}L(G_1)+(n_1+n_2)I-\frac{n_1+n_2}{n_1}J.
    \end{aligned}
\end{equation*}
Now by Lemma \ref{aL}, $M^{\dagger}= (-\frac{1}{2}L(G_1)+(n_1+n_2)I)^{-1}-\frac{1}{n_1(n_1+n_2)}J.$
Then by Lemma \ref{pinv} the $\{1\}$-inverse of
$L(G_1\dot{\vee} G_2)$ is\\

    $$\scriptsize{\begin{bmatrix}
     (-\frac{1}{2}L(G_1)+(n_1+n_2)I)^{-1}-\frac{1}{n_1(n_1+n_2)}J& (-\frac{1}{2}L(G_1)+(n_1+n_2)I)^{-1}Q_1-\frac{2}{n_1(n_1+n_2)}J   &0\\Q_1^T(-\frac{1}{2}L(G_1)+(n_1+n_2)I)^{-1}-\frac{2}{n_1(n_1+n_2)}J&\frac{1}{2}I-\frac{1}{n_1(n_1+n_2)}J+\frac{1}{4}Q_1^T(-\frac{1}{2}L(G_1)+(n_1+n_2)I)^{-1}Q_1&0\\0&0&(n_1I+L(G_2))^{-1}
    \end{bmatrix}}.$$\\
\noindent
 For any two vertices $v_i,\; v_j\in V(G_1)$,\\
        $$r_{ij}(G_1\dot{\vee} G_2)= \left(-\frac{1}{2}L(G_1)+(n_1+n_2)I\right)^{-1}_{ii}+\left(-\frac{1}{2}L(G_1)+(n_1+n_2)I\right)^{-1}_{jj}-2\left(-\frac{1}{2}L(G_1)+(n_1+n_2)I\right)^{-1}_{ij}.$$\\
For any two vertices $v_i,\; v_j\in V(G_2)$,\\
        $$r_{ij}(G_1\dot{\vee} G_2)= (n_{1}I+L(G_2))^{-1}_{ii}+(n_{1}I+L(G_2))^{-1}_{jj}-2(n_{1}I+L(G_2))^{-1}_{ij}.$$\\
For any $v_i\in V(G_1),\; v_j\in V(G_2),$\\
         $$r_{ij}(G_1\dot{\vee} G_2)= \left(-\frac{1}{2}L(G_1)+(n_1+n_2)I\right)^{-1}_{ii}+(n_{1}I+L(G_2))^{-1}_{jj}-\frac{1}{n_1(n_1+n_2)}.$$\\
From Lemma \ref{uv}, for $v_i= v_{e_{i}}\in V(C(G_1))$, where $e_i= v_pv_q\in E(G_1)$, $v_j\in V(G_1)\cup V(G_2),$\\

\begin{equation*}
    \begin{aligned}
     r_{ij}(G_1\dot{\vee} G_2)&= \frac{1}{2}\Bigl(1+r_{pj}(G_1\dot{\vee} G_2)+r_{qj}(G_1\dot{\vee} G_2)-\frac{1}{2}r_{pq}(G_1\dot{\vee} G_2)\Bigr).   
    \end{aligned}
\end{equation*}
 \\
From Lemma \ref{uv}, for $v_i= v_{e_{i}}, v_j= v_{e_{j}}\in V(C(G_1)),$ where $e_i= v_pv_q, e_j= v_sv_t\in E(G_1),$\\
 \begin{equation*}
     \begin{aligned}
         r_{ij}&= 1+\frac{1}{4}\Bigl(r_{ps}(G_1\dot{\vee} G_2)+r_{pt}(G_1\dot{\vee} G_2)+r_{qs}(G_1\dot{\vee} G_2)+ r_{qt}(G_1\dot{\vee} G_2)-r_{st}(G_1\dot{\vee} G_2)-r_{pq}(G_1\dot{\vee} G_2)\Bigr).
     \end{aligned}
 \end{equation*}

\end{proof}
\begin{theorem}\label{cej}
    Let $G_1$ be a $d$-regular graph of order $n_1$, size $m_1$ and $G_2$ be any arbitrary graph of order $n_2$. Then
    \begin{itemize}
        \item [(i)] for any two vertices $v_i,\; v_j\in V(G_1)$,\\
        \begin{equation*}
            \begin{aligned}
                r_{ij}(G_1\veebar G_2)&= (aI+bL(G_1))^{-1}_{ii}+(aI+bL(G_1))^{-1}_{jj}-2(aI+bL(G_1))^{-1}_{ij},
            \end{aligned}
        \end{equation*}
        \item[(ii)] for any two vertices $v_i,\; v_j\in V(G_2)$,\\
        $$r_{ij}(G_1\veebar G_2)= (m_{1}I+L(G_2))^{-1}_{ii}+(m_{1}I+L(G_2))^{-1}_{jj}-2(m_{1}I+L(G_2))^{-1}_{ij},$$
        \item[(iii)] for any $v_i\in V(G_1),\; v_j\in V(G_2),$\\
         $$r_{ij}(G_1\veebar G_2)= (aI+bL(G_1)^{-1}_{ii}+(m_{1}I+L(G_2))^{-1}_{jj}+\frac{n_1(n_2+2)-2d}{n_1d(n_1(n_2+2)+n_2d)},$$
        \item[(iv)] for $v_i= v_{e_{i}}, v_j= v_{e_{j}}\in V(C(G_1))$, where $e_i= v_pv_q, e_j= v_sv_t$\\
        \begin{equation*}
            \begin{aligned}
                r_{ij}(G_1\veebar G_2) &= \frac{2}{n_2+2}+\frac{1}{(n_2+2)^2}\Bigl(r_{pq}(G_1\veebar G_2)+r_{st}(G_1\veebar G_2)+4(aI+bL(G_1))^{-1}_{pq}\Bigr.\\&\Bigl.\;\;\;+4(aI+bL(G_1))^{-1}_{st}-2\left(Q(G_1)^T(aI+bL(G_1))^{-1}Q(G_1)\right)_{ij}\Bigr), 
            \end{aligned}
        \end{equation*}
        \item[(v)]  for $v_i= v_{e_{i}}\in V(C(G_1)), v_j\in V(G_1),$ where $e_i= v_pv_q$\\
        \begin{equation*}
            \begin{aligned}
                r_{ij}(G_1\veebar G_2)&= \frac{1}{(n_2+2)^2}\left(r_{pq}(G_1\veebar G_2)+4(aI+bL(G_1))^{-1}_{pq}\right)+(aI+bL(G_1))^{-1}_{jj}\\
                &\;\;\;-\frac{2}{n_2+2}\left(Q(G_1)^T(aI+bL(G_1))^{-1}\right)_{ij}+\frac{2m_1+n_2}{2m_1(n_2+2)}-\frac{n_2^2}{an_1(n_2+2)^2},
            \end{aligned}
        \end{equation*}
        \item[(vi)]for any $v_i= v_{e_{i}}\in V(C(G_1)), v_j\in V(G_2),$ where $e_i= v_pv_q$\\
        \begin{equation*}
            \begin{aligned}
                r_{ij}(G_1\veebar G_2)&= \frac{1}{(n_2+2)^2}\left(r_{pq}(G_1\veebar G_2)+4(aI+bL(G_1))^{-1}_{pq}\right)+\left(m_1I+L(G_2)\right)^{-1}_{jj}\\&\;\;\;+\frac{m_1-1}{m_1(n_2+2)}-\frac{4}{an_1(n_2+2)^2},
            \end{aligned}
        \end{equation*}

    \end{itemize}
    where $a= \frac{n_1(n_2+2)+n_2d}{n_2+2}$ and $b= -\frac{n_2+1}{n_2+2}.$
\end{theorem}

\begin{proof}
By a proper labelling of vertices in $G_1\veebar G_2$, its Laplacian matrix can be written as
$$L(G_1\veebar G_2)= \begin{bmatrix}
    n_1I-J+A(G_1)& -Q(G_1)& 0\\
    -Q(G_1)^T& (n_2+2)I& -J\\
    0& -J& m_1I+L(G_2)
\end{bmatrix}.$$
Let $L_1= n_1I-J+A(G_1), L_2= \begin{bmatrix}
    -Q(G_1)& 0
\end{bmatrix},$ and $L_3= \begin{bmatrix}
    (n_2+2)I& -J\\-J& m_1I+L(G_2)
\end{bmatrix}.$
Let \begin{equation*}
    \begin{aligned}
        M&= L_1-L_2L_3^{-1}L_2^T\\
        &= (n_1I-J+A(G_1)-\begin{bmatrix}
    -Q(G_1)& 0
\end{bmatrix}\begin{bmatrix}
    \frac{1}{n_2+2}I+\frac{n_2}{2m_1(n_2+2)}J & \frac{1}{2m_1}J\\\frac{1}{2m_1}J^T & (m_1I+L(G_2))^{-1}+\frac{1}{2m_1}J
\end{bmatrix}\begin{bmatrix}
    -Q(G_1)^T\\0
\end{bmatrix}\\&= aI+bL(G_1)-\frac{a}{n_1}J,
    \end{aligned}
\end{equation*}
where $a= \frac{n_1(n_2+2)+n_2d}{n_2+2}$ and $b= -\frac{n_2+1}{n_2+2}.$\\

\noindent Now by Lemma \ref{aL}, $M^{\dagger}= (aI+bL(G_1))^{-1}-\frac{1}{an_1}J.$\\

\noindent Then by Lemma \ref{pinv} the $\{1\}$-inverse of
$L(G_1\veebar G_2)$ is\\

    $\begin{bmatrix}M^{\dagger}&\frac{1}{n_2+2}M^{\dagger}Q(G_1)&0\\\frac{1}{n_2+2}Q(G_1)^TM^{\dagger}&\frac{1}{n_2+2}\left(I+\frac{n_2}{2m_1}J+\frac{1}{n_2+2}Q(G_1)^TM^{\dagger}Q(G_1)\right)& \frac{1}{2m_1}J\\0&\frac{1}{2m_1}J&(m_1I+L(G_2))^{-1}+\frac{1}{2m_1}J
    \end{bmatrix}.$\\
\noindent
  For any two vertices $v_i,\; v_j\in V(G_1)$,\\
        \begin{equation*}
            \begin{aligned}
                r_{ij}(G_1\veebar G_2)&= (aI+bL(G_1))^{-1}_{ii}+(aI+bL(G_1))^{-1}_{jj}-2(aI+bL(G_1))^{-1}_{ij}.
            \end{aligned}
        \end{equation*}
 For any two vertices $v_i,\; v_j\in V(G_2)$,\\
        $$r_{ij}(G_1\veebar G_2)= (m_{1}I+L(G_2))^{-1}_{ii}+(m_{1}I+L(G_2))^{-1}_{jj}-2(m_{1}I+L(G_2))^{-1}_{ij}.$$
 For any $v_i\in V(G_1),\; v_j\in V(G_2).$\\
         $$r_{ij}(G_1\veebar G_2)= (aI+bL(G_1)^{-1}_{ii}+(m_{1}I+L(G_2))^{-1}_{jj}+\frac{n_1(n_2+2)-2d}{n_1d(n_1(n_2+2)+n_2d)}.$$
 For $v_i= v_{e_{i}}, v_j= v_{e_{j}}\in V(C(G_1))$, where $e_i= v_pv_q,\; e_j= v_rv_s$\\
        \begin{equation*}
            \begin{aligned}
                r_{ij}(G_1\veebar G_2)&= \frac{2}{n_2+2}+\frac{1}{(n_2+2)^2}\Bigl(\left(Q(G_1)^T(aI+bL(G_1))^{-1}Q(G_1)\right)_{ii}\Bigr.\\
                &\;\;\;\Bigl.+\left(Q(G_1)^T(aI+bL(G_1))^{-1}Q(G_1)\right)_{jj}-2\left(Q(G_1)^T(aI+bL(G_1))^{-1}Q(G_1)\right)_{ij} \Bigr)\\
                &= \frac{2}{n_2+2}+\frac{1}{(n_2+2)^2}\Bigl(r_{pq}(G_1\veebar G_2)+r_{st}(G_1\veebar G_2)+4(aI+bL(G_1))^{-1}_{pq}+4(aI+bL(G_1))^{-1}_{st}\Bigr.\\&\Bigl.\;\;\;-2\left(Q(G_1)^T(aI+bL(G_1))^{-1}Q(G_1)\right)_{ij}\Bigr).                
            \end{aligned}
        \end{equation*}
        
\noindent For $v_i= v_{e_{i}}\in V(C(G_1)), v_j\in V(G_1),$ where $e_i= v_pv_q$\\
        \begin{equation*}
            \begin{aligned}
                r_{ij}(G_1\veebar G_2)&= \frac{1}{(n_2+2)^2}\left(Q(G_1)^T(aI+bL)^{-1}Q(G_1)\right)_{ii}+(aI+bL(G_1))^{-1}_{jj}\\
                &\;\;\;-\frac{2}{n_2+2}\left(Q(G_1)^T(aI+bL(G_1))^{-1}\right)_{ij}+\frac{2m_1+n_2}{2m_1(n_2+2)}-\frac{n_2^2}{an_1(n_2+2)^2}\\
                &= \frac{1}{(n_2+2)^2}\left(r_{pq}(G_1\veebar G_2)+4(aI+bL(G_1))^{-1}_{pq}\right)+(aI+bL(G_1))^{-1}_{jj}\\
                &\;\;\;-\frac{2}{n_2+2}\left(Q(G_1)^T(aI+bL(G_1))^{-1}\right)_{ij}+\frac{2m_1+n_2}{2m_1(n_2+2)}-\frac{n_2^2}{an_1(n_2+2)^2}.
            \end{aligned}
        \end{equation*}
    For any $v_i= v_{e_{i}}\in V(C(G_1)), v_j\in V(G_2),$ where $e_i= v_pv_q$
        \begin{equation*}
            \begin{aligned}
                r_{ij}&= \frac{1}{(n_2+2)^2}\left(Q(G_1)^T(aI+bL(G_1))^{-1}Q(G_1)\right)_{ii}+\left(m_1I+L(G_2)\right)^{-1}_{jj}+\frac{m_1-1}{m_1(n_2+2)}-\frac{4}{an_1(n_2+2)^2}\\
                &= \frac{1}{(n_2+2)^2}\left(r_{pq}(G_1\veebar G_2)+4(aI+bL(G_1))^{-1}_{pq}\right)+\left(m_1I+L(G_2)\right)^{-1}_{jj}+\frac{m_1-1}{m_1(n_2+2)}-\frac{4}{an_1(n_2+2)^2}.
            \end{aligned}
        \end{equation*}

    \noindent Here $a= \frac{n_1(n_2+2)+n_2d}{n_2+2}$ and $b= -\frac{n_2+1}{n_2+2}.$
\end{proof}
Using Lemma \ref{Kf}, Theorems \ref{cej} and \ref{cvj} we get the Kirchhoff index of central vertex (edge) join of two graphs.
\begin{theorem}\label{kfcvj}
    For $i= 1,2$, let $G_i$ be a graph of order $n_i$ and size $m_i$. Then
    \begin{equation*}
        \begin{aligned}
            Kf(G_1\dot{\vee} G_2)&= (m_1+n_1+n_2)\Biggl(tr\left((\frac{-1}{2}L(G_1)+(n_1+n_2)I)^{-1}\right)+tr\left((n_1I+L(G_2))^{-1}\right)\\&\;\;\;+\frac{1}{4}tr\left(Q(G_1)^T(\frac{-1}{2}L(G_1)+(n_1+n_2)I)^{-1}Q(G_1)\right)+\frac{m_1}{2}-\frac{m_1+n_1}{n_1(n_1+n_2)}\Biggr)\\&\;\;\;-\frac{1}{4}\alpha^T\left(\frac{-1}{2}L(G_1)+(n_1+n_2)I\right)^{-1}\alpha +\frac{m_1^2}{n_1(n_1+n_2)}-\frac{m_1}{2}-\frac{n_2}{n_1},
        \end{aligned}
    \end{equation*}
    where $\alpha= \left(d_1, d_2, \ldots, d_{n_1}\right)^T.$
\end{theorem}
\begin{theorem}\label{kfcej}
  Let $G_1$ be a $d$-regular graph of order $n_1$, size $m_1$ and $G_2$ be any arbitrary graph of order $n_2$. Then
  \begin{equation*}
      \begin{aligned}
          Kf(G_1\veebar G_2)&= (m_1+n_1+n_2)\Biggl(tr\left((aI+bL(G_1))^{-1}\right)+tr\left((m_1I+L(G_2))^{-1}\right)-\frac{4m_1}{an_1(n_2+2)^2}\\&\;\;\;+\frac{1}{(n_2+2)^2}tr\left(Q_{1}^T(aI+bL(G_1))^{-1}Q_1\right)+\frac{2m_1+n_2}{2(n_2+2)}+\frac{n_2}{2m_1}-\frac{1}{a}\Biggr)\\&\;\;\;-\frac{m_1+2n_2}{2}-\frac{n_2(n_2+2)}{2m_2}-\frac{1}{(n_2+2)^2}\alpha^TM^{\dagger}\alpha.
      \end{aligned}
  \end{equation*}
\end{theorem}
\begin{example}
    Consider $G_1\dot{\vee} G_2= C_4\dot{\vee}K_2$, then by using Theorem \ref{cvj}, the resistance matrix of $C_4\dot{\vee}K_2$ is 
$$R(C_4\dot{\vee}K_2)= \begin{bmatrix}
    0&0.45&0.4&0.45&0.61&0.61&0.81&0.81&0.37&0.37\\
    0.45&0&0.45&0.4&0.81&0.61&0.61&0.81&0.37&0.37\\
    0.4&0.45&0&0.45&0.81&0.81&0.61&0.61&0.37&0.37\\
    0.45&0.4&0.45&0&0.61&0.81&0.81&0.61&0.37&0.37\\
    0.61&0.81&0.81&0.61&0&1.1&1.2&1.1&0.75&0.75\\
    0.61&0.61&0.81&0.81&1.1&0&1.1&1.2&0.75&0.75\\
    0.81&0.61&0.61&0.81&1.2&1.1&0&1.2&0.75&0.75\\
    0.81&0.81&0.58&0.61&1.1&1.2&1.2&0&0.75&0.75\\
    0.37&0.37&0.37&0.37&0.75&0.75&0.75&0.75&0&0.33\\
     0.37&0.37&0.37&0.37&0.75&0.75&0.75&0.75&0.33&0\\
    
\end{bmatrix}$$
    Also by using Theorem \ref{kfcvj},
    $$Kf(C_4\dot{\vee}K_2)= 30.15.$$
    Similarly by using Theorems \ref{cej} and \ref{kfcej}, we get
    $$R(C_4\veebar K_2)= \begin{bmatrix}
        0&0.78&0.57&0.78&0.52&0.52&0.67&0.67&0.6&0.6\\
        0.78&0&0.78&0.57&0.67&0.52&0.52&0.67&0.6&0.6\\
        0.57&0.78&0&0.78&0.67&0.67&0.52&0.52&0.6&0.6\\
        0.78&0.57&0.78&0&0.52&0.67&0.67&0.52&0.6&0.6\\
        0.52&0.67&0.67&0.52&0&0.53&0.57&0.53&0.41&0.41\\
        0.52&0.52&0.67&0.67&0.53&0&0.53&0.57&0.41&0.41\\
        0.67&0.52&0.52&0.67&0.57&0.53&0&0.53&0.41&0.41\\
        0.6&0.6&0.6&0.6&0.41&0.41&0.41&0.41&0&0.33\\
        0.6&0.6&0.6&0.6&0.41&0.41&0.41&0.41&0.33&0
    \end{bmatrix}$$
 and   $$Kf(C_4\veebar K_2)= 25.5.$$
\end{example}

\begin{figure}[htbp]
\centering
\begin{tikzpicture}[scale=0.4]
    
\draw[ draw = black] (0,0) circle (4 cm);

  \filldraw[fill=white, draw = black] (0,4) circle (0.2 cm);
   \filldraw[fill=white, draw = black] (0,-4) circle (0.2 cm);
    \filldraw[fill=white, draw = black] (4,0) circle (0.2 cm);
     \filldraw[fill=white, draw = black] (-4,0) circle (0.2 cm);

               \filldraw[ draw = black] (2.8,2.8) circle (0.2 cm);  
     \filldraw[ draw = black] (-2.8,-2.8) circle (0.2 cm);  
        \filldraw[ draw = black] (2.8,-2.8) circle (0.2 cm);
         \filldraw[ draw = black] (-2.8,2.8) circle (0.2 cm);
          \filldraw[ draw = black] (-1.5,0) circle (0.2 cm);
           \filldraw[ draw = black] (1.5,0) circle (0.2 cm);
     \draw[-] (2.8,2.8)--(-2.8,-2.8);
     \draw[-] (-2.8,2.8)--(2.8,-2.8);
           \draw[-] (-1.5,0)--(1.5,0)--(2.8,2.8)--(-1.5,0);
           \draw[-] (1.5,0)--(-2.8,2.8)--(-1.5,0)--(-2.8,-2.8)--(1.5,0)--(2.8,-2.8)--(-1.5,0);

          \draw[ draw = black] (12,0) circle (4 cm);

               \filldraw[ draw = black] (14.8,2.8) circle (0.2 cm);  
     \filldraw[ draw = black] (9.2,-2.8) circle (0.2 cm);  
        \filldraw[ draw = black] (9.2,2.8) circle (0.2 cm);
         \filldraw[ draw = black] (14.8,-2.8) circle (0.2 cm);
          \filldraw[ draw = black] (12,2) circle (0.2 cm);
           \filldraw[ draw = black] (12,-2) circle (0.2 cm);

           \draw[-] (12,4)--(12,2);
           \draw[-] (12,-2)--(12,-4);
        \draw[-] (12,4) to [out=-135,in=135,looseness=1] (12,-2);
         \draw[-] (12,-4) to [out=45,in=-45,looseness=1] (12,2);
       \draw[-] (12,2)--(12,-2)--(16,0)--(12,2)--(8,0)--(12,-2);
       \draw[-] (14.8,2.8)--(9.2,-2.8);
       \draw[-] (14.8,-2.8)--(9.2,2.8);
        \filldraw[fill=white, draw = black] (16,0) circle (0.2 cm);
   \filldraw[fill=white, draw = black] (8,0) circle (0.2 cm);
    \filldraw[fill=white, draw = black] (12,4) circle (0.2 cm);
     \filldraw[fill=white, draw = black] (12,-4) circle (0.2 cm);
      
\end{tikzpicture}
\caption{$C_4\dot{\vee}K_2$ and $C_4\veebar K_2$.}
 \end{figure}

\section{Conclusion}
    In this article, the concepts of resistance distance in the central graph and the central vertex (edge) join of two graphs are explored. Furthermore, using the results, the Kirchhoff index of the graphs is determined.
\section{Declarations}
 On behalf of all authors, the corresponding author states that there is no conflict of interest.

 \bibliography{cvej}
 \bibliographystyle{plain}
\end{document}